\documentclass[12pt]{amsart}

\usepackage{tikz-cd}

\usepackage[english]{babel}
\usepackage[utf8x]{inputenc}
\usepackage[T1]{fontenc}
\usepackage{amsmath}
\usepackage{amstext}
\usepackage{amsfonts}
\usepackage{amssymb}
\usepackage{amsthm}
\usepackage{amsrefs}
\usepackage{mathtools}
\usepackage{dsfont}
\usepackage{color}
\usepackage{graphicx}
\usepackage[colorinlistoftodos]{todonotes}
\usepackage[colorlinks=true, allcolors=blue]{hyperref}
\usepackage{float}
\usepackage[colorinlistoftodos]{todonotes}
\usepackage{theoremref}
\usepackage{comment}
\usepackage{tikz-cd}
\allowdisplaybreaks
\newtheorem{theorem}{Theorem}
\newtheorem{lemma}[theorem]{Lemma}
\newtheorem{prop}[theorem]{Proposition}

\newtheorem*{cor*}{Corollary}
\newtheorem*{thm*}{Theorem}
\newtheorem*{lem*}{Lemma}

\newtheorem*{prop*}{Proposition}
\theoremstyle{definition}
\newtheorem{definition}[theorem]{Definition}

\newtheorem{example}[theorem]{Example}
\newtheorem*{defn*}{Definition}

\theoremstyle{remark}

\newcommand{\Ad}{\operatorname{Ad}}

\newcommand{\bE}{E}

\newcommand{\aut}{\operatorname{{Aut}}}

\newcommand{\ev}{\text{ev}}
\newcommand{\res}{\text{res}}

\title{Powers averaging for actions on $C(X)$-algebras}
\author{Tattwamasi Amrutam}
\address{Ben Gurion University of the Negev.
	Department of Mathematics.
	Be'er Sheva, 8410501, Israel.}
\email{tattwama@post.bgu.ac.il}
\author{Ilan Hirshberg}
\address{Ben Gurion University of the Negev.
	Department of Mathematics.
	Be'er Sheva, 8410501, Israel.}
\email{ilan@bgu.ac.il}
\author{Apurva Seth}
\address{Ben Gurion University of the Negev.
	Department of Mathematics.
	Be'er Sheva, 8410501, Israel.}
\email{apurva@post.bgu.ac.il}
\date{\today}

\thanks{We thank the Fields Institute, in which most of the work on this paper was done, for its hospitality}
\begin{document}
\maketitle
\begin{abstract}
Given a unital $C(X)$-algebra $A$, a discrete group $\Gamma$, and an action $\alpha \colon\Gamma\to \aut(A)$ which leaves  $C(X)$ invariant and such that $C(X)\rtimes_{\alpha,r} \Gamma$ is simple, and a $2$-cocycle $\omega$, we obtain a bijective correspondence between maximal $\Gamma$-invariant ideals of $A$ and maximal ideals in $A\rtimes_{\alpha,\omega,r} \Gamma$. In particular, $A\rtimes_{\alpha,\omega,r} \Gamma$ is simple if and only if $A$ has no $\Gamma$-invariant ideals.
\end{abstract}

Let $X$ be a compact Hausdorff space and $A$ be an unital $C(X)$-algebra. We consider $C(X)$ a distinguished subalgebra of the center $Z(A)$.
 Let $\Gamma$ be a discrete group and $\alpha \colon \Gamma \to \aut(A)$ be an action. We say that $\alpha$ is \emph{compatible} with an action on $C(X)$ if $C(X)$ is invariant under $\alpha$. Our main result in this paper is the following. It generalizes \cite[Theorem~7.1]{BKKO} and \cite[Theorem 1.1]{bryder2018reduced}, which are special cases of Theorem~\ref{gensimplicity} when $X = \{ \text{pt} \}$. By a \emph{maximal ideal} in the statement below, we mean a maximal proper ideal.

\begin{theorem}
\thlabel{gensimplicity}
Let $X$ be a compact Hausdorff space and $A$ be an unital $C(X)$-algebra. Let $\Gamma$ be a discrete group, and let $\alpha \colon \Gamma \to \aut(A)$ be an action that is compatible with an action on $C(X)$, also denoted by $\alpha$.  Suppose $C(X) \rtimes_{\alpha,r} \Gamma$ is simple. Let $\omega$ be a $2$-cocycle for the action.
Let $I_{\max}(A,\Gamma)$ denote the set of maximal $\Gamma$-invariant ideals in $A$, and let $I_{\max}(A \rtimes_{\alpha,\omega,r} \Gamma)$ denote the set of maximal ideals in $A \rtimes_{\alpha,\omega,r} \Gamma$. For an invariant ideal $J$ in $A$, we denote by $\tilde{\alpha}$ the induced action on $A/J$. Then the maps
\[
\iota \colon I_{\max}(A,\Gamma) \to I_{\max}(A \rtimes_{\alpha,\omega,r} \Gamma)
\]
given by 
\[
\iota (J) = \ker \left ( A \rtimes_{\alpha,\omega,r} \Gamma \to ( A/J ) \rtimes_{\tilde{\alpha},\omega,r} \Gamma \right )
\]
and 
\[
\res \colon I_{\max}(A \rtimes_{\alpha,\omega,r} \Gamma)  \to I_{\max}(A,\Gamma)  \text{ given by }
\res (J) = J \cap A
\]
define a bijective correspondence between those sets of maximal ideals. In particular, if there are no $\Gamma$-invariant ideals, then $A \rtimes_{\alpha,\omega,r} \Gamma$ is simple.
\end{theorem}

The assumption that $C(X) \rtimes_{\alpha,r} \Gamma$ implies that the action of $\Gamma$ on $X$ is minimal. However, we do not assume that the action on $X$ is topologically free. For example, if $X = \{ \text{pt} \}$, so the action on $X$ is trivial, this simply means that $\Gamma$ is $C^*$-simple, which is the assumption in  \cite[Theorem~7.1]{BKKO} and \cite[Theorem 1.1]{bryder2018reduced}. 

\begin{example}
	Let $X$ be a minimal $\Gamma$-space such that $\Gamma\curvearrowright X$ is topologically free. Denote by $\alpha$ the action on $C(X)$. It follows from the result on page 122 of \cite{archbold-spielberg} that $C(X)\rtimes_{\alpha,r} \Gamma$ is simple. Suppose now that $D$ is an unital $C^*$-algebra, suppose $\beta \colon \Gamma \to \aut(D)$ is an action with no $\Gamma$-invariant ideals. Let $A=D \otimes C(X)$. It follows from \thref{gensimplicity} that $A\rtimes_{\beta \otimes \alpha , r} \Gamma$ is simple. 
\end{example}

We begin by recalling some definitions and results. For a compact Hausdorff space $X$, we use $\delta_x$ to denote the point probability measure at $x$, and by $\ev_x$ the evaluation state at $x$. We use both notations as sometimes it is more convenient to talk about states and sometimes about measures. When we have an action of $\Gamma$ on a space $X$, and the action is understood, we may write $(s,x) \mapsto s \cdot x$ for the action with $s \in \Gamma$ and $x \in X$. We write $\Gamma x = \{ s \cdot x : g \in \Gamma\}$. Likewise, if $\nu$ is a measure on $X$ and $s \in \Gamma$, we write $s_* \nu$ for the push-forward measure, and $\Gamma \nu = \{ s_*  \nu : s \in \Gamma\}$.
\subsection*{Generalized Boundaries and Contractible Measures}

The notion of a generalized boundary was introduced independently in \cite{Naghavi} and \cite{Kawabe} as a generalization of the notion of a $\Gamma$-boundary (which Furstenberg~\cite{Furstenberg1973} introduced and has been used extensively recently, see \cite{KK, BKKO, HartKal}).
Given a $\Gamma$-space $X$ and a probability measure $\nu\in\text{Prob}(X)$, we say that $\nu$ is contractible if
\[\left\{\delta_x : x\in X \right\}\subset \overline{\Gamma\nu}^{\text{w}^*}.\] 
Notice that any Dirac measure is contractible in a minimal $\Gamma$-space.
\begin{definition}[Generalized boundary]
Let $\Gamma$ be a discrete group. Let $X$ and $Y$ be compact Hausdorff spaces endowed with an action of $\Gamma$. Suppose that the action on $X$ is minimal. Let $\pi: Y\to X$ be a factor map, that is, a continuous surjective $\Gamma$-equivariant map. We say that $Y$ is a $(\Gamma,X)$-boundary if 
\begin{enumerate}
    \item The action of $\Gamma$ on $Y$ is minimal. 
    \item For any probability measure $\nu\in \text{Prob}(Y)$, if $\pi_*\nu$ is contractible then so $\nu$. 
\end{enumerate}
\end{definition}
Given a minimal $\Gamma$-space $X$, there exists a maximal $(\Gamma,X)$-boundary, which is the spectrum of the $\Gamma$-injective envelope of $C(X)$ (see for example \cite[Section~3]{Kawabe}). We denote this space by $\partial_F(\Gamma, X)$. It follows from \cite[Theorem~3.4]{Kawabe} that if the action on $X$ is minimal, then $C(X)\rtimes_{\alpha,r}\Gamma$ is simple if and only if $\Gamma\curvearrowright\partial_F(\Gamma, X)$ is topologically free (that is, for any $s \in \Gamma \smallsetminus \{e\}$, the set of points $\partial_F(\Gamma, X)$ not fixed by $s$ is dense). By  \cite[Proposition~3.3]{Kawabe}, since $\partial_F(\Gamma, X)$ is a Stonean space, an action on $\partial_F(\Gamma, X)$ is topologically free if and only if it is free.
\subsection*{Tensor Products of $C(X)$-Algebras over $C(X)$}
We briefly recall the construction of tensor products of $C(X)$-algebras over $C(X)$. We refer the reader to \cite{blanchardtensor} for further details. 
 Let $A$ be a $C^*$-algebra, and $X$ be a compact Hausdorff space. A $C(X)$-algebra consists of a $C^*$-algebra $A$ along with a unital $*$-homomorphism $\iota \colon C(X) \to Z(M(A))$.  ($Z(M(A))$ denotes the center of the multiplier algebra of $A$.) We shall assume throughout that $\iota$ is injective. We suppress the notation for $\iota$ and think of $C(X)$ as a distinguished subalgebra of $Z(M(A))$. We focus in this paper only on unital $C^*$-algebras, so $A = M(A)$ in our case.

Let $A$ and $B$ be two $C(X)$-algebras, and let $I$ be the involutive ideal of the algebraic tensor product $A\otimes_{\text{alg}}B$ generated by the elements $(fa)\otimes b- a\otimes (fb)$, for $f\in C(X)$, $a\in A$, and $b\in B$. Then the quotient algebra $A\otimes_{\text{alg}}B/I$ is an involutive algebra over $\mathbb{C}$. 

If $A$ and $B$ are unital $C(X)$-algebras, and one of them is abelian, then there is a unique $C^*$-norm on $A\otimes_{\text{alg}}B/I$ (\cite[Lemma 2.7]{blanchardtensor}). Note that in this case, $I$ and $J$ (defined as in \cite[Definition~2.1]{blanchardtensor}) coincide (see for example, \cite[Proposition~3.1]{blanchardtensor}). This is the only case we consider in this paper, so we denote it by $A\otimes_{C(X)}B$.

\subsection*{Crossed Products} Recall that a (scalar-valued) $2$-cocycle for a discrete group $\Gamma$ is a function $\omega \colon \Gamma \times \Gamma \to \mathbb{T} \subset \mathbb{C}$ which satisfies, for any $r,s,t \in \Gamma$,
\[
\omega(r,s)\omega(rs,t) = \omega(r,st)\omega(s,t) 
\]
and
\[
\omega(e,t) = \omega(t,e) = 1 \, .
\]
(One can consider more general unitary-valued cocycles for an action, provided those unitaries are in the normalizer of $C(X)$, but we do not go into the added generality here, as we do not have interesting enough examples to justify it.)
Given an unital $C^*$-algebra, an action $\alpha \colon \Gamma \to \aut(A)$ and a $2$-cocycle $\omega$, the twisted reduced crossed product $A \rtimes_{\alpha,\omega,r} \Gamma$ is the subalgebra of the algebra of adjointable operators on the Hilbert $C^*$-module $l^2(\Gamma) \otimes A$ (with $A$ given the structure of a trivial $C^*$-module over $A$) generated by the operators $\{ a \lambda_s: a \in A, s \in \Gamma\}$ defined on elements of the form $\delta_t \otimes b$ (for $t \in \Gamma$ and $b \in A$) by
\[
a \lambda_s \cdot \delta_t \otimes b = \omega((st)^{-1},s)\delta_{st} \otimes \alpha_{(st)^{-1}}(a)b \, .
\]
We view $A$ as embedded in $A \rtimes_{\alpha,\omega,r} \Gamma$ in the canonical way. That is, the unitaries $\{\lambda_s : s \in \Gamma\}$ satisfy $\lambda_s\lambda_t = \omega(s,t) \lambda_{st}$ and $\lambda_s a \lambda_s^* = \alpha_s(a)$ for any $s,t \in \Gamma$ and for any $a \in A$. 
We denote by $E: A \rtimes_{\alpha,r} \Gamma \to A$ the canonical faithful conditional expectation which is given by $E(a\lambda_s) = \delta_{s,e} a$ for any $a\in A$ and for any $s \in \Gamma$. Notice that the map $s \mapsto \Ad_{\lambda_s}$ is an action of $\Gamma$ on $A \rtimes_{\alpha,\omega,r}\Gamma$  by automorphisms, and the map $E$ is $\Gamma$-equivariant with respect to this inner action of $\Gamma$.

Note that if $J$ is an ideal of $A \rtimes_{\alpha,\omega,r} \Gamma$, then $J \cap A$ is a $\Gamma$-invariant ideal in $A$.

\noindent
\subsection*{Generalized Probability Measures}
The notion of a generalized probability measure was introduced in \cite{AU}. We recall the definition for the sake of completion.
\begin{definition}
\thlabel{definition_generalized_probability_measure}
	Let $X$ be a compact Hausdorff space, $\Gamma$ be a discrete group, and $\alpha \colon \Gamma \to \aut(C(X))$ be an action. A \emph{generalized $(C(X),\Gamma,\alpha)$-probability measure} is a function $f \colon \Gamma \to C(X)_+$ such that $\sum_{s \in \Gamma} f(s)^2 = 1$ (where the sum converges uniformly). 
	
	We denote the set of all such generalized measures by $P(\Gamma, C(X),\alpha)$. Moreover, $P_f(\Gamma, C(X),\alpha)$ denotes the collection of generalized probability measures with finite support. 
\end{definition}
It is sometimes convenient to consider sums over a (typically finite) subset $I \subset \Gamma$, in which 
	case, it is understood that the elements corresponding to $s \in \Gamma \smallsetminus I$ are zero. 
	
    Given a $C(X)$-algebra $A$ along with an action $\alpha \colon \Gamma \to \aut(A)$ leaving $C(X)$ invariant, and given any $\mu \in P(\Gamma,C(X),\alpha)$, we define a unital and completely positive map $\Phi_{\mu} \colon A \to A$ by
    $$ \Phi_{\mu} (a) = \sum_{s \in \Gamma} f(s) \alpha_s( a )  f(s) \, . $$
    The map $\Phi_{\mu}$ induces a dual map on the state space $\Phi_{\mu}^* \colon S(A) \to S(A)$.
    
   We extend the definition of $\Phi_{\mu}$ to a completely positive map 
   \[
   \Phi_{\mu} \colon A \rtimes_{\alpha,\omega,r} \Gamma \to A \rtimes_{\alpha,\omega,r} \Gamma 
   \]
    by setting 
   \[
   \Phi_{\mu} (a) = \sum_{s \in \Gamma} f(s) \lambda_s  a \lambda_s^*  f(s) \, .
   \] 
	Likewise, $\Phi_{\mu}^*$ extends to an affine map
	\[
	\Phi_{\mu}^* \colon S( A \rtimes_{\alpha,\omega,r} \Gamma ) \to  S( A \rtimes_{\alpha,\omega,r} \Gamma )  \, .
	\]
	
To prove the simplicity result, we need an averaging result, which was also used in \cite{Haagerup, HartKal, AU}, that shows that probability measures can be contracted to Dirac measures using $C(X)$-convex-combination of group elements. 
\begin{lemma}\cite[Lemma~3.6]{AU}
\thlabel{lemma_generalized_contractibility}
    Let $X$ be a compact Hausdorff space and $\Gamma$ be a discrete group. Suppose $\alpha \colon \Gamma \to C(X)$ is an action that induces a minimal action on $X$. For any $x \in X$ there is a net $(\mu_\lambda) \subseteq P_f(\Gamma,C(X),\alpha)$ with the property that for any $\nu \in S(C(X))$, we have $\Phi_{\mu_\lambda}^*(\nu) \to \ev_x$.
\begin{proof}
Let $U$ be an open neighborhood of $x$.  For any $x_0\in X$, the orbit $\Gamma x_0$ is dense in $X$ since the action is assumed to be minimal. Therefore, we can find an element $t_{x_0}\in G$ such that $t_{x_0}x_0\in U$, or equivalently, $x_0 \in t_{x_0}^{-1} U$. Since $X=\cup_{t \in \Gamma }t^{-1}U$, by compactness, we can find a finite set $t_1,t_2,\ldots , t_n \in \Gamma$ such that  $t_{1} U, \dots, t_{n}U$ are a cover of $X$. Now let $\{g_{i}\}_{i=1,2,\ldots n}$ be a partition of unity subordinate to $\left\{t_{i}U\right\}_i$. Define $f \colon \{ t_{1} , t_{2} , \ldots ,  t_{n} \} \to C(X)_+$ by 
$ f (t_{i} ) = \sqrt{g_{i}} $.
Given any $\nu \in P(X)$, it follows that $\Phi^*_{\mu_U} ( \nu )$ is a measure with support contained in the closure of $U$. The net $\left\{\Phi^*_{\mu_U} \nu\right\}$, indexed by open neighborhoods of $x$ ordered by reverse inclusion, converges to $\ev_x$ in the weak$^*$-topology, as required.
\end{proof}    
\end{lemma}
\begin{lemma}
\thlabel{singularity}
Let $X$ be compact Hausdorff space, let $A$ be an unital $C(X)$-algebra, let $\Gamma$ be a discrete group, and let $\alpha \colon \Gamma \to \aut(A)$ be an action which leaves $C(X)$ invariant, and such that the induced action on $X$ is minimal. Let $\omega$ be a $2$-cocycle for $\Gamma$. Fix $x \in X$. 
Let $\eta\in S(A \rtimes_{\alpha,\omega,r} \Gamma)$ be such that $\eta|_{C(X)}=\ev_x$. Then, for every $s\not\in\Gamma_x$ and for any $a \in  ( A \rtimes_{\alpha,\omega,r} \Gamma  )\cap C(X)' $ we have $\eta(a \lambda_s)=0$.
\begin{proof}
Since $\eta|_{C(X)}=\ev_x$, the subalgebra $C(X)$ in the multiplicative domain of $\eta$. Let $s\in \Gamma$ be such that $sx\ne x$. Using Uryhson's lemma, choose $f \in C(X)$ such that $f(x) = 1$ and $f(sx) = 0$. Fix an element $a \in ( A \rtimes_{\alpha,\omega,r} \Gamma ) \cap C(X)' $. We have
\begin{align*}\eta(a\lambda_s)&=f(x) \eta(a\lambda_s)\\&= \eta(f a\lambda_s)\\&= \eta(af\lambda_s
)\\&= \eta(a\lambda_s) f(sx)\\&=0. \end{align*}
\end{proof}    
\end{lemma}

\begin{prop}
\thlabel{onemeasure}
Let $X$ be compact Hausdorff space, let $A$ be an unital $C(X)$-algebra, let $\Gamma$ be a discrete group,  let $\alpha \colon \Gamma \to \aut(A)$ be an action which leaves $C(X)$ invariant.
Suppose that $C(X) \rtimes_{\alpha,r} \Gamma$ is simple. Let $\omega$ be a $2$-cocycle for $\Gamma$. Then, for any state $\varphi\in S(A\rtimes_{\alpha,\omega,r} \Gamma)$ there exists a state $\psi \in S(A)$ such that
$$ \psi \circ E\in\overline{\{\Phi_{\mu}^* ( \varphi ) : \mu \in P_f(C(X),\Gamma,\alpha)\}}^{\text{w}^{*}} \, . $$
\begin{proof}
Let $I_{\Gamma}(C(X)) = C(\partial_F(\Gamma,X))$ be the $\Gamma$-injective envelope of $C(X)$. Denote by $\beta$ the canonical extension of the action of $\Gamma$ to $ C(\partial_F(\Gamma,X))$.  Let $\tilde{\alpha} = \alpha \otimes_{C(X)} \beta$ on $A \otimes_{C(X)}  C(\partial_F(\Gamma,X))$. 
Extend the state $\varphi$ to a state $\tilde{\varphi}$ on $\left(A\otimes_{C(X)} C(\partial_F(\Gamma,X))\right) \rtimes_{\tilde{\alpha},\omega,r} \Gamma$. Because $C(X) \rtimes_{\alpha,r} \Gamma$ is simple, the action on $X$ is minimal.  
Therefore, using \thref{lemma_generalized_contractibility}, we can find a net $(\mu_\lambda) \subseteq P_f(\Gamma,C(X),\alpha)$ such that $\Phi_{\mu_\lambda}^* ( \varphi|_{C(X)} ) \to \ev_x$ for some $x \in X$. Passing to a subnet if needed, we may assume that there exists 
\[
 \tilde{\psi} \in S \left ( \left(A\otimes_{C(X)} C(\partial_F(\Gamma, X))\right) \rtimes_{\tilde{\alpha},\omega,r} \Gamma \right )
 \]
  such that $\Phi_{\mu_\lambda}^* ( \tilde{\varphi} ) \to \tilde{\psi}$. Notice that the state $\tilde{\psi}$ satisfies $\tilde{\psi}|_{C(X)} = \ev_x$. 
    
Therefore, $\tilde{\psi}|_{C(\partial_F(\Gamma,X))}$ is contractible as well.
 Hence, using \cite[Theorem~A]{Naghavi}, there is a net $(s_i)_{i \in I} \subset \Gamma$ with $\tilde{\psi}|_{C(\partial_F(\Gamma,X))} \circ \beta_{s_i} \to \ev_y$ for some $y \in \partial_F(\Gamma,X)$. Again passing to a subnet if needed, we can find a state $\eta \in \overline{\{ \Phi_{\mu}^* ( \tilde{\varphi} ): \mu \in P_f(\Gamma, C(X),\alpha)\}}^{\text{w}^{*}}$ with the property that $\eta|_{C(\partial_F(\Gamma, X))} = \ev_y$. 
    
We claim that $\eta|_{A \rtimes_{\alpha,\omega,r} \Gamma}$ satisfies the requirement in the statement. By \cite[Theorem~3.4 and Proposition 3.3]{Kawabe}, because $C(X) \rtimes_{\alpha,r} \Gamma$ is simple, the action of $G$ on $\partial_F(\Gamma,X)$ is free. 
 Observe that $C(\partial_F(\Gamma,X))$ is in the commutant of $A$ in $A \otimes_{C(X)} C(\partial_F(\Gamma,X))$. It then follows from \thref{singularity} that for any $a \in A$, we have $\eta(a \lambda_s) = 0$ whenever $s \neq e$. Thus $\eta|_{A \rtimes_{\alpha,r} \Gamma} = \eta|_{A \rtimes_{\alpha,r} \Gamma} \circ E$, as required.
\end{proof}
\end{prop}
The following generalizes \cite[Lemma~4.1]{bryder2018reduced}. 
\begin{lemma}
\thlabel{properideal}
Let $X$ be a compact Hausdorff space, let $A$ be an unital $C(X)$-algebra, let $\Gamma$ be a discrete group, and let $\alpha \colon \Gamma \to \aut(A)$ be an action which leaves $C(X)$ invariant. Suppose $C(X) \rtimes_{\alpha,r} \Gamma$ is simple. Let $\beta$ be the canonical extension of $\alpha|_{C(X)}$ to $C(\partial_F(\Gamma,X))$, and let $\tilde{\alpha} = \alpha \otimes_{C(X)} \beta$ be the induced action on 
$A \otimes_{C(X)} C(\partial_F(\Gamma,X))$. 
Let $\omega$ be a $2$-cocycle for $\Gamma$. If $I$ is proper ideal in $A\rtimes_{\alpha,\omega,r} \Gamma$, then the ideal $J$ generated by $I$ inside $(A\otimes_{C(X)}C(\partial_F(\Gamma,X))\rtimes_{\tilde{\alpha},\omega,r} \Gamma$ is proper. 
\begin{proof}
Let $\varphi$ be a state on $A\rtimes_{\alpha,\omega,r} \Gamma$ such that $\varphi|_{I}=0$. Let $\tilde{\varphi}$ be an extension of this state to $A\otimes_{C(X)}C(\partial_F(\Gamma,X)) \rtimes_{\tilde{\alpha},\omega,r}  \Gamma$.
Following a similar argument to that of proof of \thref{onemeasure}, there exists a state $\psi$ on $A\otimes_{C(X)}C(\partial_F(\Gamma,X)) \rtimes_{\tilde{\alpha},\omega,r}  \Gamma$, such that $\psi|_{C(\partial_F(\Gamma,X))}=\ev_y$ for some $y \in \partial_F(\Gamma,X)$ and such that
$$ \psi \in\overline{\{\Phi_{\mu}^* ( \tilde{\varphi} ) : \mu \in P_f(C(X),\Gamma,\alpha)\}}^{\text{w}^{*}} \, . $$
Note that $C(\partial_F(\Gamma,X))$ is in the multiplicative domain of $\psi$, and $\psi|_{I}=0$. 

Now, for any element $x\in I$, for any $a,b \in A$, for any $f,g \in C(\partial_F(\Gamma,X))$ and for any $s,t\in \Gamma$, because $\psi (axb) = 0$, we have 
\begin{align*}
&\psi\left((a\otimes_{C(X)} f )\lambda_{s} x(b \otimes g)\lambda_t\right)\\&=f_1(y)\psi\left(a \lambda_{s} )xb\lambda_t \right)g (t \cdot y)\\&=0 \, .
\end{align*}
Therefore we have $\psi|_J=0$. Hence, $J$ is a proper ideal. 
\end{proof}
\end{lemma}
Let $X$, $A$, $\Gamma$ and $\omega$ be as above.
Given a $\Gamma$-invariant ideal $I\triangleleft A$, let $I\rtimes_{\alpha,\omega,r}\Gamma$ be the ideal in $A\rtimes_{\alpha,\omega,r}\Gamma$ generated by $I$. We write $\overline{\alpha}$ for the induced action of $\Gamma$ on $A/I$. 
The surjection $\pi^I \colon A \to A/I$ induces a  surjective $*$-homomorphism $\tilde{\pi}^I \colon A\rtimes_{\alpha,\omega,r} \Gamma \to A/I \rtimes_{\overline{\alpha},\omega,r} \Gamma$. 
We do not necessarily have $I\rtimes_{\alpha,\omega,r}\Gamma = \ker (\tilde{\pi}^I)$. This happens when the group $\Gamma$ is exact; see \cite{kirchberg_wassermann_exact_groups,exel_exact}. 
We write $E_I$, $E_A$ and $E_{A/I}$ for the canonical expectation maps 
\[
E_I \colon I\rtimes_{\alpha,\omega,r}\Gamma  \to I \, , 
\]
\[
E_A \colon A\rtimes_{\alpha,\omega,r}\Gamma  \to A 
\]
and 
\[
E_{A/I} \colon A\rtimes_{\overline{\alpha},\omega,r}\Gamma  \to A/I \, .
\]
 We have a commuting diagram
\begin{equation}
\label{eq:firstcommutativediagram}
\begin{tikzcd}
&
I\rtimes_{\alpha,\omega,r}\Gamma  \arrow[r] \arrow[d,"E_I"'] &
A\rtimes_{\alpha,\omega,r}\Gamma  \arrow[r,"\tilde{\pi}^I"] \arrow[d,"E_A"'] &
 A/I\rtimes_{\overline{\alpha},\omega,r}\Gamma \arrow[d,"E_{A/I}"'] &
\\
0 \arrow[r] &
I \arrow[r] &
A \arrow[r,"\pi"] &
A/I \arrow[r] &
0
\end{tikzcd}
\end{equation}
It follows from diagram~(\ref{eq:firstcommutativediagram}) that  
\begin{equation}
\label{eq:equalityofidealsandintersection}
\ker(\tilde{\pi}^I) \cap A=I.    
\end{equation}
Now, given an ideal $J\triangleleft(A\otimes_{C(X)}C(\partial_F(\Gamma,X))\rtimes_{\tilde{\alpha},\omega,r} \Gamma$, write $J_1 = J \cap A$ and $J_2 = J \cap A\otimes_{C(X)}C(\partial_F(\Gamma,X))$. 
We have the following commutative diagram of $*$-homomorphisms. 
\[ \begin{tikzcd}
A\rtimes_{\alpha,\omega,r} \Gamma \arrow{r}{}  \arrow[d,"\tilde{\pi}^{J_1}"']  & (A\otimes_{C(X)}C(\partial_F(\Gamma,X))\rtimes_{\tilde{\alpha},\omega,r}\Gamma  \arrow[d,"\tilde{\pi}^{J_2}"'] \\%
( A/J_1) \rtimes_{\overline{\alpha},\omega,r} \Gamma \arrow{r}{}&  \left [
A\otimes_{C(X)} C(\partial_F(\Gamma,X)) / J_2 ) 
\right ] \rtimes_{\overline{\tilde{\alpha} },\omega,r}\Gamma
\end{tikzcd}
\]
The horizontal arrows are injective. Therefore,
\begin{equation}
    \label{eq:equalityafterlifting}
\ker(\tilde{\pi}^{J_1})  = \ker(\tilde{\pi}^{J_2}) \cap A\rtimes_{\alpha,\omega,r} \Gamma \, .
\end{equation}

The following is a generalization of \cite[Lemma~7.3]{BKKO}. 
\begin{lemma}
\thlabel{propercontainment}
Let $X$ be a compact Hausdorff space, let $A$ be an unital $C(X)$-algebra, let $\Gamma$ be a discrete group, and let $\alpha \colon \Gamma \to \aut(A)$ be an action which leaves $C(X)$ invariant. Suppose $C(X) \rtimes_{\alpha,r} \Gamma$ is simple. Let $\beta$ be the canonical extension of $\alpha|_{C(X)}$ to $C(\partial_F(\Gamma,X))$, and let $\tilde{\alpha} = \alpha \otimes_{C(X)} \beta$ be the induced action on 
$A \otimes_{C(X)} C(\partial_F(\Gamma,X))$. 
Let $\omega$ be a $2$-cocycle for $\Gamma$. 
 Let $J$ be an ideal in $\left(A\otimes_{C(X)}C(\partial_F(\Gamma,X))\right) \rtimes_{\tilde{\alpha},\omega,r} \Gamma$. Let $J_2$, $\tilde{\pi}^{J_1}$ and $\tilde{\pi}^{J_2}$ be as in the paragraph before the statement of the lemma.   Then, 
\[ 
J_2 \rtimes_{\tilde{\alpha},\omega,r} \Gamma
\subseteq J \subseteq 
\ker (\tilde{\pi}^{J_2})
\, .
\]
\begin{proof}
	The inclusion $J_2 \rtimes_{\tilde{\alpha},\omega,r} \Gamma \subseteq J$ is immediate. We show the other inclusion. 
	
	For $x \in X$, we denote by $A_x$ the fiber of $A$ over $x$, that is, $A_x = A / C_0(X \setminus \{x\})A$. 
	Let us denote the factor map by $p: \partial_F(\Gamma, X)\to X$. Notice that for any $y \in \partial_F(\Gamma,X)$,
	 $\text{id}\otimes_{C(X)}\ev_y$ maps  $A\otimes_{C(X)}C(\partial_F(\Gamma,X))$ into $A_{p(y)}$.
Let 
\[
J_2=J\cap\left(A\otimes_{C(X)}C(\partial_F(\Gamma,X))\right) \, .
\]
Fix $y \in \partial_F(\Gamma,X)$. Consider the map 
\[
\pi_y=\text{id}\otimes_{C(X)}\ev_y: A\otimes_{C(X)}C(\partial_F(\Gamma,X))\xrightarrow{} A_{p(y)}
\, .
\]
 Let $J_2^y=\text{id}\otimes_{C(X)}\delta_y(J_2)$. 
We note that $J_2^y$ is an ideal of $A_{p(y)}$. Therefore, the u.c.p map 
\[
\Theta_y: A\otimes_{C(X)}C(\partial_F(\Gamma,X))/J_2\xrightarrow{}A_{p(y)}/J_2^y
\]
given by
$\Theta_y ( b+J_2 ) =  \pi_y(b)+J_2^y$ is well-defined. Pick now a Hilbert space $K$ and an embedding $j \colon A_{p(y)}/J_2^y \to B(K)$. 
Consider the composition 
\begin{align*}J+A\otimes_{C(X)}C(\partial_F(\Gamma,X))&\xrightarrow{}\frac{J+A\otimes_{C(X)}C(\partial_F(\Gamma,X))}{J}\\&=\frac{A\otimes_{C(X)}C(\partial_F(\Gamma,X))}{J_2}\\&\xrightarrow{\Theta_y}\frac{A_{p(y)}}{J_2^y}\xrightarrow{j} B (K),\end{align*}
By Arveson's extension theorem, we can choose a u.c.p map 
\[
\Psi_y: \left(A\otimes_{C(X)}C(\partial_F(\Gamma, X))\right) \rtimes_{\tilde{\alpha},\omega,r} \Gamma\to B(K)
\]
which extends this composition. Let 
\[
 E :\left({A\otimes_{C(X)}C(\partial_F(\Gamma, X)}\right) \rtimes_{\tilde{\alpha},\omega,r}\Gamma\to {A\otimes_{C(X)}C(\partial_F(\Gamma, X)}
 \] 
 be the canonical conditional expectation.
 We claim that $\Psi_y=\Psi_y\circ E $. Observe that ${A\otimes_{C(X)}C(\partial_F(\Gamma, X)}$ is in the multiplicative domain of $\Psi_y$. Moreover, $\Psi_y(f)=f(y)$ for all $f\in C(\partial_F(\Gamma, X))$. Therefore, to show that $\Psi_y=\Psi_y\circ E $, it is enough to show that $\Psi_y( \lambda_s )=0$ for all $s\in\Gamma\setminus\{e\}$. Since $C(X)\rtimes_{\alpha,r}\Gamma$ is simple, $\Gamma\curvearrowright\partial_F(\Gamma, X)$ is free, and hence, $\Gamma_y=\{e\}$. Since $sy\ne y$, we can find $f\in C(\partial_F(\Gamma, X))$ such that $f(y)=1$ and $f(sy)=0$. Thus, \begin{align*}\Psi_y( \lambda_s )&=f(y)\Psi_y( \lambda_s )\\&=\Psi_y(f \lambda_s )\\&=\Psi_y( \lambda_s \alpha_{s^{-1}} ( f ))\\&=\Psi_y( \lambda_s )f(sy)\\&=0 \end{align*}
Consequently, we see that $\Psi_y(J)=0=\Psi_y( E (J))$ for all $y\in \partial_F(\Gamma, X)$. Since $ E (J)\subset A\otimes_{C(X)}C(\partial_F(\Gamma, X)$, 
it follows that $ E (J)\subset J_2$. This is equivalent to saying that $J \subset \ker (\tilde{\pi}^{J_2})$. 
\end{proof}
\end{lemma}

We can now prove  \thref{gensimplicity}.
\begin{proof}[Proof of  \thref{gensimplicity}]
Adopt the notation from the discussion before \thref{propercontainment}. Let $I$ be a maximal $\Gamma$-invariant ideal of $A$. Let $J$ be a proper ideal in $A \rtimes_{\alpha,\omega,r} \Gamma$ such that $\ker (\tilde{\pi}^I) \subseteq J$. Let $\hat{J}$ be the ideal generated by $J$ in $\left(A\otimes_{C(X)}C(\partial_F(\Gamma,X))\right) \rtimes_{\tilde{\alpha},\omega,r} \Gamma$. Let 
\[
L = \hat{J} \cap A\otimes_{C(X)}C(\partial_F(\Gamma,X)) \, .
\]
By \thref{propercontainment}, we have $\hat{J} \subseteq \ker {\tilde{\pi}^L}$. 
Using equation~\eqref{eq:equalityafterlifting}, we see that
\begin{equation}
	\label{eq:onesideinclusion}
	J \subseteq 
	( \hat{J} \cap A ) \rtimes_{\alpha,\omega,r} \Gamma \subseteq
	\left ( \ker ( \tilde{\pi}^L ) \cap A \right )  \rtimes_{\alpha,\omega,r} \Gamma  =
	\ker ( \tilde{\pi}^{\hat{J} \cap A} ) \,. 
	\end{equation}
By applying equation~\eqref{eq:equalityofidealsandintersection} to $I$ and to $\hat{J} \cap A$, we see that
\[
I \subseteq J\cap A \subseteq \hat{J}\cap A \, .
\]
Since $J$ is a proper ideal, from \thref{properideal}, it follows that $\hat{J}$ is proper. Therefore, from the maximality of $I$, we obtain  $I=\hat{J} \cap A$. Combining this along with equation~\eqref{eq:onesideinclusion}, we see that $J\subset \ker ( \tilde{\pi}^{\hat{J} \cap A} ) = \ker ( \tilde{\pi}^{I} ) = \iota (I)$. This shows that $\iota(I)$ is a maximal ideal inside $A \rtimes_{\alpha,\omega,r} \Gamma$.

For the other direction, let $J$ be a maximal ideal inside $A \rtimes_{\alpha,\omega,r} \Gamma$. We must show that $I =J\cap A$ is a maximal $\Gamma$-invariant ideal on $A$. Let $\hat{J}$ be the ideal generated by $J$ inside $(A\otimes_{C(X)}C(\partial_F(\Gamma,X)) \rtimes_{\tilde{\alpha},\omega,r}  \Gamma$. 
 Let $L = \hat{J} \cap A\otimes_{C(X)}C(\partial_F(\Gamma,X))$. 
Using  \thref{propercontainment}, we see that 
$\hat{J} 
\subset 
\ker ( \tilde{\pi}^{L})$. 

Combining this with equation~\eqref{eq:equalityafterlifting}, we obtain that
\[
J \subseteq \hat{J} \cap  ( A\rtimes_{\alpha,\omega,r}\Gamma)
\subseteq 
\ker ( \tilde{\pi}^{L}) \cap  ( A\rtimes_{\alpha,\omega,r}\Gamma ) =
\ker ( \tilde{\pi}^{\hat{J} \cap A} )  \, .
\]
Since $J$ is a proper ideal in $A \rtimes_{\alpha,\omega,r} \Gamma$, using \thref{properideal}, we see that $\hat{J} \cap A$ must be a proper ideal inside $A$. By the maximality of $J$, it follows that $J=\ker ( \tilde{\pi}^{\hat{J} \cap A} ) $. Using equation~\eqref{eq:equalityofidealsandintersection}, we see that 
\[
I =J\cap A=
\ker ( \tilde{\pi}^{\hat{J} \cap A} ) \cap A=\hat{J} \cap A \, .
\]
Consequently, we see that $J=\ker (\tilde{\pi}^I)$. 
Now, if $\Tilde{I}$ is a $\Gamma$-invariant proper ideal of $A$ containing $I$, then we obtain that 
$J = \ker (\tilde{\pi}^I) \subseteq \ker \left ( \tilde{\pi}^{\tilde{I}} \right )$.
Because $ \ker \left ( \tilde{\pi}^{\tilde{I}} \right )$ is a proper ideal in $A\rtimes_{\alpha,\omega,r}\Gamma$ which contains $J$, and $J$ is maximal, it follows that
$J =  \ker \left ( \tilde{\pi}^{\tilde{I}} \right )$. Consequently,
\[
I = J \cap A=  \ker \left ( \tilde{\pi}^{\tilde{I}} \right ) \cap A= \tilde{I} \, ,\]
as required.
\end{proof}

\bibliographystyle{amsalpha}
\bibliography{Simplicity_bib.bib}

\end{document}